\documentclass[12pt]{amsart}
\usepackage[utf8]{inputenc}
\usepackage{amssymb}
\usepackage{amsmath}
\usepackage{amsthm}
\usepackage{graphicx}
\usepackage{amscd}
\usepackage{epic}
\usepackage[alphabetic]{amsrefs}
\usepackage{xypic}
\usepackage[colorlinks, citecolor=blue]{hyperref}
\usepackage[left=1.5in,right=1.5in,top=1.5in,bottom=1.5in]{geometry}
\usepackage{soul}
\usepackage[pagewise]{lineno}
\usepackage{comment}
\usepackage{enumitem}
\usepackage{mathrsfs}

\theoremstyle{plain}
\newtheorem{theorem}{Theorem}[section]

\newtheorem{lemma}[theorem]{Lemma}

\newtheorem{proposition}[theorem]{Proposition}
\newtheorem{definition-lemma}[theorem]{Definition-Lemma}
\theoremstyle{remark}

\theoremstyle{definition}

\newcommand{\supp}[0]{\operatorname{Supp}}

\def\>{\geq}
\newcommand{\mbQ}{\mathbb{Q}}
\newcommand{\mbR}{\mathbb{R}}

\newcommand{\mult}{{\rm mult}}
\newcommand{\Supp}{{\rm Supp}}

\newcommand{\mbC}{\mathbb{C}}

\def\injective{\hookrightarrow}

\def\mbP{\mathbb{P}}
\def\>{\geq}

\def\mult{\operatorname{mult}}

\def\Ex{\operatorname{Ex}}
\def\dim{\operatorname{dim}}

\def\Chow{\operatorname{Chow}}

\theoremstyle{definition}
\newtheorem{definition}[theorem]{Definition}
\theoremstyle{definition}

\numberwithin{equation}{section}

\theoremstyle{remark}

\author{Roktim Mascharak}

\title{On Log Canonical MMP for K\"ahler 3-Fold}

\begin{document}
\maketitle
\begin{abstract}
    In this article we show that the Log Minimal Model Program holds for $\mbQ$-factorial lc pair $(X,\Delta)$ with $X$ being a compact K\"ahler $3$-fold having only klt singularities.
\end{abstract}
\tableofcontents
\section{Introduction}
In algebraic geometry one of the fundamental question has been to classify the  algebraic varieties up-to birational equivalence. One of the most important tool in birational classification of complex projective variety is the minimal model program. It was fully established in dimension $3$ in $80$'s and  $90$'s, and recently extended to many cases in arbitrary dimension including the case of varieties of log general type \cite{BCHM}.\\
Unfortunately we don't know much about minimal model program for compact K\"ahler varieties. In recent years by the works of \cite{HP16,CHP,HP15} existence and termination of the MMP for compact K\"ahler varieties with terminal singularity was proved , and due to \cite{DH22}, \cite{DO22}, we know informations about log minimal model program for compact K\"ahler $3$-folds upto pairs with dlt singularities. In this paper we try to extend the minimal model program for the most general class of MMP singularities, i.e. Log canonical (lc) singularities. One of the problems with varieties with log canonical singularities are that lc singularity is not rational singularity. But for K\"ahler varieties if we do not assume that it has rational singularity, we loose many valuable properties so that we can not even set up the convex geometrical premises to run the minimal model program. Also to run the MMP and do divisorial contraction we will need to assume that $X$ is klt. So in this article we take $X$ to be a compact K\"ahler $3$-fold having rational singularities and $X$ is klt. Next we state our main theorems. \\
We start by proving the cone theorem for K\"ahler $3$-folds with lc singularities. More precisely the following.
 \begin{theorem}\label{C1}
Let $(X,B)$ be a compact $K\ddot{a}ler$ 3-fold lc pair with $B$ effective and $X$ having rational singularity. Then there exist a countable collection of rational curves $\lbrace C_i \rbrace$ on $X$ such that 
\begin{enumerate}
\item $\overline{NA}(X)=\overline{NA}(X)_{(K_X+B)\geqslant0} + \sum \mathbb{R}^+\cdot [C_i]$
\item $-6\leqslant (K_X+B)\cdot C_i<0$
\item for any K\"ahler class $\omega$, $(K_X+B+\omega)\cdot C_i \leqslant 0$ for all but finitely many $i$ 
\item The rays $\lbrace \mathbb{R}^+\cdot [C_i] \rbrace$ do not accumulate in $(K_X+B)_{<0}$
\end{enumerate}
\end{theorem}
 To proof this theorem we go to a dlt model and then push forward the dlt cone theorem. Note that to prove this theorem we don't need to assume $X$ has only klt singularities.\\ 
We prove existence of divisorial contraction (only for the case when the whole divisor is getting contracted to a point) and flipping contractions for compact K\"ahler $3$-fold lc pairs $(X,\Delta)$ with rational singularities, with $K_X+\Delta$ pseudo-effective as well as non pseudo-effective. But to prove the divisiorial contraction in case when divisor is getting contracted to a curve, we will need to assume that $X$ is of klt type, i.e. $(X,0)$ is klt. We will elaborate more on this in section $5$. Also we have shown existence of flips for compact K\"ahler $3$-fold lc pairs using the same method as in the case of complex projective $3$-fold lc pairs. More precisely we prove the following statement.
\begin{theorem}\label{Cm}
    Let $X$ be a compact K\"ahler klt $3$-fold. Let $\Delta$ be a boundary divisor such that $(X,\Delta)$ is a lc pair. Let $R$ be an $K_X+\Delta$-negative extremal ray of $\overline{NA}(X)$. Then the contraction of $R$, $cont_R:X\rightarrow Z$ exists.
\end{theorem}
 Next we prove that if the starting variety $X$ has rational singularities then each step of the MMP preserves rational singularities.
 \begin{theorem}\label{T4}
 
 Let $(X,\Delta)$ be compact K\"ahler $3$-fold lc pair with $X$ having rational singularities. Then the following holds
\begin{enumerate}

\item let $f:X\rightarrow Y$ be a divisorial contraction of $(K_X+\Delta)$ negative extremal ray. Then $Y$ also has rational singularity.
\item Let $(X,\Delta)$ be a lc pair. $f:X\rightarrow Y$ be a $K_X+\Delta $-flipping contraction and $f^+:X^+\rightarrow Y$ be the flip. If $X$ has rational singularity then $X^+$ has rational singularity.

\end{enumerate}
\end{theorem}
 To prove the above theorems we will follow Fujino's treatment as in the the article \cite{FV22}.\\
 In Section $6$ we prove the termination of sequence of flips in this set up. which is as follows.
 \begin{theorem}\label{T1}
Let $(X,\Delta)$ be a compact K\"ahler $3$-fold lc pair with $\mbQ$ factorial rational singularity. Then every sequence of $K_X+\Delta$ flips terminate.
	
\end{theorem}

In the next section we develop the MMP for compact K\"ahler $3$-fold lc pair $(X,\Delta)$ with $\mbQ$-factorial rational singularities and $K_X+\Delta$ not pseudo-effective, i.e. state and proof cone and contraction theorems accordingly. In the last section we proof the main theorem in the above set-up which is as follows.
\begin{theorem}\label{FT}
Let $(X,B)$ be a lc pair where $X$ is a $\mbQ$-factorial compact K\"ahler $3$-fold, with $X$ having klt singularity. If $K_X+B$ is not pseudo effective then there exists a finite sequence of flips and divisorial contractions\[
	\Phi:X\dashrightarrow X_1\dashrightarrow...\dashrightarrow X_n.
\]
and a Mori fiber space $\phi:X_n\rightarrow S$, i.e. a morphism such that $-(K_{X_n}+\Phi_* B)$ is $\phi$-ample and $\rho(X_n/S)=1$	

\end{theorem}
 Combining all the above theorems we have the following theorem
 \begin{theorem}\label{MMP}
     Let $(X,\Delta)$ be a compact K\"ahler $3$-fold lc pair with $X$ having klt singularity. Then we can run the $K_X+\Delta$-MMP and the following holds
     \begin{enumerate}
         \item If $K_X+\Delta$ is pseudo-effective then the sequence of $K_X+\Delta$-flips and divisorial contraction terminates with a minimal model $(X',\Delta')$, i.e. $K_{X'}+\Delta'$ is nef.
         \item If $K_X+\Delta$ is not pseudo-effective, then the $K_X+\Delta$-MMP terminates with a Mori fiber space.
     \end{enumerate}
 \end{theorem}
\textbf{Acknowledgement:} The author would like to thank Professor Omprokash Das for his guidance and many insightful discussions.
\section{Preliminaries}

An analytic variety or simply a variety is an irreducible reduced complex space. In this article by a pair $(X,B)$ we will always mean that $X$ is normal, $B\geqslant 0$ is an effective $\mbQ$-divisor such that $K_X+B$ is $\mbQ$-Cartier. 
\begin{definition}\cite[Definition $2.2$]{HP16}
	An analytic variety $X$ is called K\"ahler if there exists a positive closed $(1,1)$ form $\omega\in \mathcal{A}_{\mbR}^{1,1}(X)$ such that the following holds: for every point $x\in X$ there exists an open neighbourhood $x\in U$ and a closed embedding $i_U:U\injective V$ into an open set $V\in \mbC^N$, and a smooth strictly plurisubharmonic function $f:V\rightarrow \mbR$ such that $\omega|_{U\cap X_{sm}}=(i\partial\overline{\partial}f )|_{U\cap X_{sm}}$. Here $X_{sm}$ is the smooth locus of $X$.

\end{definition}

Now to do MMP for cases of K\"ahler variety we will need some definitions to set up the perfect counterpart for cone of curves and other machinery in this scenario. In the following we collect some important definition regarding that. More details on these can be found in \cite{HP16}, \cite{CHP} and the references therein.
\begin{definition}
	\begin{enumerate}
		\item A compact variety $X$ is said to belong to \textit{Fujiki's class} $\mathcal{C}$ if 
	 $X$ is bimeromorphic to a compact K\"ahler manifold is satisfied.\\
\\
\item On a normal compact analytic variety $X$ we replace the use of N\'eron-Severi group $NS(X)_{\mbR}$ by $H^{1,1}_{BC}(X)$, the Bott-Chern cohomology of real closed $(1,1)$ forms with local potentials or equivalently, the closed bi-degree $(1,1)$ currents with local potentials. See \cite[Definition $3.1$]{HP16} for more details. So we define \[
	N^1(X) = H^{1,1}_{BC}(X).
\]
\item If $X$ is in Fujiki's class $\mathcal{C}$, then from \cite[Eqn. 3]{HP16} we know that $N^1(X)=H^{1,1}_{BC}(X)\subset H^2(X,\mbR)$. In particular, the intersection product can be defined in $N^1(X)$ via the cup product of $H^2(X,\mbR)$. 
\item Let $X$ be a normal compact variety in Fujiki's class $\mathcal{C}$. We define $N_1(X)$ to be the vector space of real closed currents of bi-dimension $(1,1)$ modulo the following equivalence relation: $T_1\cong T_2$ if and only if $T_1(\eta)=T_2(\eta)$ for all real closed $(1,1)$ forms $\eta$ with local potentials. 
\item For the definitions of closed cone of currents $\overline{NA}(X)$, pseudoeffective currents, and nef class we reffer to \cite{HP16}.
\item The nef cone $Nef(X)\subset N^1(X)$ is the cone generated by nef cohomology classes. Let $\mathcal{K}$ be the open cone in $N^1(X)$ generated by the classes of K\"ahler forms. Nef cone is closure of $\mathcal{K}$( see \cite[Remark $3.12$]{HP16}).
\item We say that a variety $X$ is $\mbQ$-factorial if every weil divisor $D\subset X$, there is a positive intiger $k>0$ such that $kD$ is a Cartier divisor, and for the canonical sheaf $\omega_X$, there is a positive intiger $m>0$ such that $(\omega_X^{\otimes m})^{**}$ is a line bundle. It is well known that if $X$ is a $\mbQ$-factorial $3$-fold and $X\dashrightarrow X'$ is a flip or divisorial contraction, then $X'$ is also $\mbQ$-factorial.

	\end{enumerate}
	
\end{definition}

\section{Cone theorem}	
	
Let $X$ be $\mathbb{Q}$-factorial Kahler 3-Fold, $\Delta$ be an effective-divisor, such that $(X,\Delta)$ has dlt singularity. Then we have the following theorem 
\begin{theorem}\cite[Theorem $2.17$,Theorem $4.6$]{DH22}
Let $(X,B)$ be a $\mathbb{Q}$-factorial compact $K\ddot{a}ler$ 3-fold dlt pair. Then there exists a countable collection of rational curves $\lbrace C_i \rbrace _{i\in I}$ such that $0<-(K_X+B)\cdot C_i \leqslant 6$ and $$\overline{NA}(X)= \overline{NA}(X)_{(K_X+B)\geqslant 0} +\sum \mathbb{R}^+ \cdot [C_i]  $$ Moreover if $\omega$ is a $K\ddot{a}hler$ class, then there are only finitely many extremal rays $R_i=\mathbb{R}^+ .[C_i]$ satisfying $(K_X+B+\omega)\cdot C_i <0$ for $i\in I$.
\end{theorem}
Now we will try to extend the above theorem for $K\ddot{a}ler$ 3-fold lc pair. But before that we will need another theorem about the faces of $\overline{NA}(X)$. 
\begin{lemma}\cite[Lemma 3.1]{JW17}\label{cg}
Let $f:V\rightarrow W$ be a surjective linear map of finite dimensional vector spaces. Suppose $C_V\subset V$ and $C_W \subset W$ are closed convex cones of maximal dimension and $H\subset W$ is a linear subspace of codimension 1 assume that $f(C_V)=C_W$ and $C_W\cap H \subset \partial C_W$ then $f^{-1}H\cap C_V=f^{-1} (H\cap C_W)\cap C_V$ and $f^{-1}H\cap C_V \subset \partial C_V$.
\end{lemma}

Now we come to the proof of cone theorem for K\"ahler $3$-fold with lc singularities.

\begin{proof}[Proof of Theorem \ref{C1}]
$(X,B)$ is given to be a log canonical pair. Take $f:Y\rightarrow X$ be a log resolution of $(X,B)$. so consider $B'$ to be strict transform of $B$. Then we have the log equation $K_Y+B'=f^*(K_X+B)+\sum a_i E_i$. As $(X,B)$ is log canonical pair, $a_i\geqslant -1$. Let $F= \sum E_i$ where $E_i$ are all exceptional divisor of $f$ such that $a_i=-1$ in the above log equation. Now adding $F$ to both side of the equation we get that the divisor $B'+F=B_Y$ and $\sum a_i E_i +F$ are both effective . Then running relative MMP w.r.t. the dlt pair $(Y,B_Y)$ over $X$, by special termination and negativity lemma we have a pair $(\bar Y,\bar B _Y)$ such that this a $\mathbb{Q}$-factorial dlt pair with $K_{\bar Y} + \bar B_Y =\bar f ^*(K_X+ B)$. Now replacing $(Y,B_Y)$ by $(\bar Y,\bar B _Y)$ consider the linear map $f_* : N_1(Y)\rightarrow N_1(X)$, since if we can take a curve $C\subseteq f_*^{-1} C'$ , where $C'\in N_1(X)$ ,then $f_* C=C'$, hence the map $f_*$ is surjective. Note that we can do this since $f$ is a projective morphism. So we have $$f_*(\overline{NA}(Y)) = \overline{NA}(X). $$ By the $\mathbb{Q}$-factorial dlt cone theorem, there is a countable collection of rational curves $C_i ^Y$ on $Y$ such that $$\overline{NA}(Y)= \overline{NA}(Y)_{K_Y+B_Y \geqslant 0} + \sum \mathbb{R}_{\geqslant 0} . [C_i ^Y]$$. Let $C_i$ be the countable collection of the curve on $X$ given by $f_* C_i ^Y$ with reduced structure. We claim that these curves satisfy (1) in theorem. Suppose instead  $$\overline{NA}(X)\neq \overline{NA}(X)_{(K_X+B)\geqslant 0} +\sum \mathbb{R}^+ .[C_i]  $$. Then there is some cohomology class $\alpha$ which is positive on right hand side of the equation but non positive on $\overline{NA}(X)$. Let $\omega$ be an K\"ahler class and $\lambda$= inf$\lbrace t: \alpha+t\omega$ is nef$ \rbrace$ then $\alpha+\lambda \omega$ is nef but not ample. So by Kleimann's criterion it takes value 0 somewhere on $\overline{NA}(X) \setminus 0$. Replacing $\alpha$ by $alpha+\lambda \omega$ we may assume that $\alpha$ is non negative on $\overline{NA}(X)\setminus 0$. But since this is non K\"ahler hence $\alpha=0$ intersects $\overline{NA}(X)\setminus 0$ non trivially, so it cuts out some face of $\overline{NA}(X)$. By the Lemma \ref{cg}  we know $F_Y= f_* ^{-1} F\cap \overline{NA}(Y)$ is some non empty extremal face of $\overline{NA}(Y)$. Which is $K_Y+B_Y$-negative away for a lower dimensional face $f_*^{-1}(0)$. Now by \cite[Theorem 1.7]{DH22} we know that for any such class $\alpha$ nef and $\alpha-(K_Y+B_Y)$ nef we have a projective surjective morphism $\Phi:Y\rightarrow Z$ over $X$ with $\alpha=\Phi^*(\alpha_Z)$, where $\alpha_Z$ is a $K\ddot{a}ler$ class. So this is a analogue of base-point free theorem. Hence any such face $F_Y$ contains a negative extremal ray which is generated by one of the $C_i^Y$. But then $\alpha_{=0}$ contains one of the $C_i$ which is a contradiction.\\
To prove the inequality $$-6\leqslant(K_X+B)\cdot C_i<0$$ We first notice that $0<\frac{(K_X+B)\cdot C_i}{(K_X+B)\cdot f_* C_i^Y}\leqslant 1$, this follows from the definition of $C_i$ and. Now from $\mathbb{Q}$-factorial dlt cone theorem we know that $$-6\leqslant(K_Y+B_Y)\cdot C_i^Y<0$$ Hence by projection formula and and the above inequality we have our desired inequality. \\
Next we show that the rays $R_i= \mathbb{R}_{\geqslant 0} [C_i] $ do not accumulate in $(K_X+B)_{<0}$. Suppose otherwise. So there is some sequence $\lbrace R_i\rbrace$ converging to a $K_X+B$ negative ray $R$. Let $R_I^Y$ be  extremal rays in $\overline{NA}(Y)$ satisfying $f_* R_i^Y=R_i$. By compactness of the unit ball in $\overline{NA}(Y)$, some sub-sequence of $R_i^Y$ must converge to a ray $R^Y$. By continuity of $f_*$ we must have $f_*(R_Y)=R$. Now we know that $K_{Y} + B_Y = f ^*(K_X+ B)$ and $(K_X+B).R<0$ so by projection formula we have that $(K_Y+B_Y).R_i^Y<0$. But this contradicts the cone theorem for K\"ahler $\mathbb{Q}$ factorial dlt pairs.\\
Finally let $\omega$ be an K\"ahler class on $X$. Suppose there are infinitely many $C_i$ with $(K_X+B+\omega)\cdot C_i<0$. By compactness, some subsequence of the corresponding $R_i$s converge to a ray $R$. This must satisfy $(K_X+B+\omega).R\leqslant 0$. But since $R\subset \overline{NA}(X)$ so $(K_X+B)\cdot R\leqslant -\omega\cdot R <0$ as $\omega\cdot R>0$, which contradicts the part (4) of this theorem.

\end{proof}
\section{Contraction Theorems}
In this section we will proof the contraction theorems for $(X,\Delta)$ lc $3$-fold pair with $X$ compact K\"ahler and $K_X+\Delta$ pseudo-effective. First we will talk about divisorial contractions. We a priori have the cone theorem for $(X,\Delta)$. We say that a $(K_X+\Delta)$ negative extremal ray is small if every curve $[C]\in R$ has the property that $\dim_{mC} \Chow(X)=0$ for all natural number $m$. So in a small ray no curves deformes non trivially, and since Chow space has countably many components there can be at most countably many distinct curves in a small ray, hence  the locus covered by all the curves in a small ray can not be a surface.\\
A $(K_X+\Delta)$ negative extremal ray is called divisorial if it is not small. Using the similar method as in proof of \cite[Lemma 5.6]{HP16} we can show that for any curve $[C]$  in a divisorial ray $\dim_C \Chow(X)>0$.

\begin{theorem}\label{Cont}
Let $X$ be a compact K\"ahler $3$-fold with at-most rational singularity, $\Delta$ be a boundary divisor such that $(X,\Delta)$ is a lc pair with $K_X+\Delta$ pseudo effective. Let $R$ be an $K_X+\Delta$ negative extremal ray of divisorial type. Then there exists a bimeromorphic morphism $f:X\rightarrow Y$ which contracts the locus which is covered by curves in $R$.	
	
\end{theorem}
\begin{proof}
$R$ is an extremal ray of $\overline{NA}(X)$, hence there exists a $(1,1)$ nef class $\alpha$ such that $\alpha^{\bot}\cap \overline{NA}(X)=R$. For a K\"ahler form $\omega$ we have $\overline{NA}(X)=\overline{NA}(X)_{K_X+\Delta+\omega\geqslant 0}+\sum_{finite} \mbR^+[C_i]$. We can choose $\omega$ in such way so that $R=\mbR^+[C_{i0}]$ is also a $K_X+\Delta+\omega$ negative ray. Hence in particular $\overline{NA}(X)_{K_X+\Delta+\omega\geqslant 0}+\sum_{i\neq i0} \mbR^+[C_i]$ is also a closed cone, so with proper scaling of $\alpha$ we have that $\alpha-(K_X+\Delta)$ is a K\"ahler class. So $\alpha$ is also a big class. Now we claim that $K_X+\Delta \sim_{\mbQ} \sum \lambda_j S_j$ where $\lambda_j\geqslant 0$ and $S_j$ are irreducible surfaces. To see this first we go to the dlt model of the lc pair $(X,\Delta)$ which exists by \cite[Theorem 5.2]{DO22}, let that be $(X',\Delta')$ and this a crepant blow-up so we have $h:X'\rightarrow X$ such that $K_{X'}+\Delta'=h^*(K_X+\Delta)$. Now we run the MMP for the dlt pair $(X',\Delta')$ and it terminates with $(X'',\Delta'')$ with $K_{X''}+\Delta''$ nef, since $K_{X'}+\Delta'$ was pseudo-effective. This MMP terminates by \cite[Theorem $1.1$]{DH22}. Now let $\phi:X'\dashrightarrow X''$ be the birational morphism, we take $W$ to be the normalization of the graph of $\phi$. So we have the two maps $p:W\rightarrow X'$ and $q:W\rightarrow X''$, and we have that $p^*(K_{X'}+\Delta')=q^*(K_{X''}+\Delta'')+E$ with $E$ effective, by negativity lemma. We have the non vanishing theorem for K\"ahler $3$-fold dlt pair by \cite[Theorem $1.1$]{DO22}, so $K_{X''}+\Delta''$ is $\mbQ$-effective, so is $p^*(K_{X'}+\Delta')=q^*(K_{X''}+\Delta'')+E$. Hence by projection formula $K_{X'}+\Delta'$ is $\mbQ$-effective and so is $K_X+\Delta$ as $h$ is a crepant blowup. Now $(K_X+\Delta)\cdot R<0$ implies that $\sum \lambda_j S_j \cdot R<0$ which means that $S_{j_0}.R<0$ for some $j_0$. In particular $S_{j_0}\cdot [C]<0$ for any $[C]\in R$. Let $\lbrace C_t \rbrace_{t\in T}$ is a deformation family of $C$. Consider $\cup_{t\in T}C_t=A$. Since all the $C_t$ are numerically equivalent we have that $A\subset S_{j_0}$. $A$ can not be a curve since $\dim_C \Chow(X)>0$, which means $C$ deforms non-trivially. Hence $\overline{A}$ is a surface which is contained in irreducible surface $S_{j_0}$. So closure of $A$ is $S_{j_0}$ itself. Hence $S_{j_0}$ is the unique surface which contains all the deformation family of any curve class in $R$. Now $X$ is $\mbQ$-factorial hence there exists a $m$ such that $mS_{j_0}$ is a Cartier divisor. Now let $\mu:S\rightarrow S_{j_0}$ be normalization of $S_{j_0}$. So $\mu^*(\alpha|_{S_{j_0}})$ is a nef and big class. So $S\rightarrow Z$ be the nef reduction of $\mu^*(\alpha|_{S_{j_0}})$. Since $\alpha$ intersects a covering family of curves of $S_{j_0}$ trivially, $\mu^*(\alpha|_{S_{j_0}})$ Intersects a family of curves which coveres a dense open subset of $S$. Hence $Z$ can have dimension $0$ or $1$. We first deal with the case when $\dim Z=0$ hence $Z$ is a point. So we have a map from $S_{j_0}\rightarrow Z$. We will show that this morphism extends to give a contraction of $S_{j_0}$. It is sufficient to show that $\mathcal{O}_{S_{j_0}}(-mS_{j_0})$ is an ample line bundle by \cite[Proposition $7.4$]{HP16}. Now $-mS_{j_0}$ is strictly positive on $R$ and $\alpha$ is strictly positive on $\overline{NA}(X)_{K_X+\Delta+\omega}+\sum_{i\neq i0} \mbR^+[C_i]$. So there exists an $\epsilon>0$ such that $\alpha-\epsilon mS_{j_0}$ is strictly positive on $\overline{NA}(X)$, which implies $\alpha-\epsilon mS_{j_0}=\omega$ for a K\"ahler class $\omega$.
Now $\alpha|_{S_{j_0}}$ is numerically trivial . So $-mS_{j_0}|_{S_{j_0}}=\frac{1}{\epsilon} \omega|_{S_{j_0}}$ is an apmle Cartier divisor. Hence  $\mathcal{O}_{S_{j_0}}(-mS_{j_0})$ is an ample line bundle. So by \cite[Proposition $7.4$]{HP16} we have our desired contraction $f:X\rightarrow Y$.\\
Now we come to the case when $dimZ=1$, i.e. $Z$ is an curve. Now we will follow the same approach as \cite[proof of Theorem $1.5$]{DH22}. So we will go to a log resolution $\mu:X'\rightarrow X$ and set $\Delta'=\mu_*^-1(\Delta+(1-b)S_{j_0})+\Ex(\mu)$ where $b= \mult_{S_{j_0}}\Delta$. So $(X',\Delta')$ is a dlt pair. 

Let $\alpha '=\mu^*(\alpha)$. Now just as in the proof of \cite[Theorem 1.5]{DH22} we run a $\alpha'$-trivial $K_X'+\Delta'$ MMP and we get \[
	\phi_i :X'=X_0\dashrightarrow X_1'\dashrightarrow X_2'-\dashrightarrow X_i'
\] 
This MMP exists and terminates in the exact same way as described it the above mentioned paper. Suppose in terminates at $X_n'$, then the main issue is even though we can show that $\phi_{n*}(\sum E_i+S')=0$ for all $E_i$ whose discrepancy is strictly greater than $=-1$, and $S'$ being strict transform of $S_{j_0}$. But in this way we are unable to contract those exceptional divisor $F_i$ which are over lc centers. So if we assume $X$ to be klt then for $\epsilon<1$, $(X,(1-\epsilon)\Delta)$ is klt. Now for small enough $\epsilon$, $R$ is also $K_X+(1-\epsilon)\Delta$-negative extremal ray, and hence we can contract it.

\end{proof}

Now we consider the case when $R$ is a $K_X+\Delta$ negative extremal ray which is not of the type divisorial type, i.e. $R=\mbR^+[C]$ with $\dim_{mC} \Chow(X)=0$ for all $m$ natural number.
So we will show the existence of flipping contractions. 
\begin{theorem}\label{C'1}
Let $(X,\Delta)$ be an $3$-fold lc pair with $X$ compact K\"ahler. $R=\mbR^+[C]$ be a $K_X+\Delta$ negative extremal ray with $[C]$ very rigid. Then there exists a bimeromorphic morphism $f:X\rightarrow Y$ which contracts exactly the locus covered by curve classes inside $R$, and is an isomorphism outside that locus.	
	
\end{theorem}

\begin{proof}
	If $X$ is not uniruled then existence of such $f$ readily follows from \cite[Theorem $4.14$]{CHP}. If $X$ is uniruled the main obstruction is to show that- if $\alpha$ is the supporting nef class of $R$, and $S$ is the locus covered by curves in $R$ then $(\alpha|_S)^2=0$. But the same proof as \cite[Proposition 4.12]{DH22} works for lc pair $(X,\Delta)$ as proof does not use any property of dlt pair. Now the existence of contraction $f$ again follows from \cite[Theorem 4.14]{CHP}. 
	
\end{proof}
Next we proof the existence of log canonical flips in this scenario.
\begin{theorem}
Let $(X,\Delta)$ be a log canonical pair where $X$ is a compact K\"ahler $3$-fold having rational singularities, $f:X\rightarrow Z$ be a $K_X+\Delta$ negative small contaction. Then the flip of $f$ exists.	
	
\end{theorem}
\begin{proof}
The proof is exactly same as \cite[Theorem $8.1$]{AST}. The main two ingredients of this proof is relative basepoint free theorem and relative K-V vanishing theorem, whose analogues in K\"ahler variety set up are respectively \cite[Theorem $4.8$]{Nak} and \cite[Theorem 2.16]{DH22}. 
\end{proof}

\section{Rational Singularity}
Now we will suppose that $X$ has rational singularity otherwise we will loose many valuable information like dulaity of closure of K\"ahler cone of currents and nef cone. So assuming existence of divisorial contractions and flip we will show that the target variety also has rational singularity.

\begin{proof}[Proof of Theorem \ref{T4}](1)
To prove the theorem first we will use \textit{Reid-Fukuda type vanishing theorem} to show that $R^qf_*\mathcal{O}_X=0$. Note that $f$ is a projective morphism. Now $(X,\Delta)$ is an lc pair, we claim that for small $\epsilon>0$ the pair $(X,(1-\epsilon)\Delta)$ does not have any log canonical center along the support of $\Delta$. Suppose it has a lc center along support of $\Delta$ then let $E$ be the irreducible exceptional divisor with $a(E,X,(1-\epsilon)\Delta)=-1$. Then $a(E,X,\Delta)=a(E,X,(1-\epsilon)\Delta+\epsilon \Delta)=-1-c$, where $c$ is the coeffecient of $E$ in pull back of $\epsilon\Delta$. Which contradicts the fact that $(X,\Delta)$ is a lc pair. So $(X,(1-\epsilon)\Delta)$ is klt along support of $\Delta$. Now we can choose $\epsilon$ small enough so that $f$ is also a contraction of $K_X+(1-\epsilon)\Delta$ negative extremal ray. As the image of exceptional divisors are curves or points hence there are no lc centre of $(X,(1-\epsilon)\Delta)$ which is a surface. $-(K_X+(1-\epsilon)\Delta)$ is $f$ ample, $C$ be an lc cnter which is a curve. So $f(C)$ is either a point or a curve. If $f(C)$ is a point then $C$ is a $(K_X+(1-\epsilon)\Delta)$ negative curve hence $-(K_X+(1-\epsilon)\Delta)|_C$ is ample. If $f(C)$ is a curve then then fibers of $f|_C$ are finitely many point, hence restriction of $-(K_X+(1-\epsilon)\Delta)$ is again ample, so it is $f$ ample over $f(C)$. Hence by \cite[Theorem $5.8$]{FCC} applying to the pair $(X,(1-\epsilon)\Delta)$ we have that $R^if_*\mathcal{O}_X=0$ for all $i>0$. Now we will use Leray spectral sequence to proof that $Y$ also has rational singularity. Let $Z$ be an common resolution of $X$ and $Y$. So let the corresponding maps be $g:Z\rightarrow X$ and $h:Z\rightarrow Y$. It is enough to show that $h$ is a rational resolution, apriory by Zariski's main theorem we have that $h_*\mathcal{O}_Z=\mathcal{O}_Y$, so we just need to show that $R^ih_*\mathcal{O}_Z=0$ for all $i>0$. We have that $f_*\circ g_*$ and $h_*$ are same functor from category of sheaves of coherent sheaves over $Z$ to category of coherent sheaves over $Y$. So by Leray Spectral sequence we have that $R^if_*\circ R^jg_* \Rightarrow R^{i+j}h_*$. We already know that $R^ig_*\mathcal{O}_Z=0$. Hence in the $E_2$ page of the spectral sequence $R^if_*\circ R^jg_*$ applied on the sheaf $\mathcal{O}_Z$, we see that only non zero entry is in $(0,0)$ position which is $\mathcal{O}_Y$. So the spectral sequence terminates at the $E_2$ page and hence we have $R^0h_*\mathcal{O}_Z=\mathcal{O}_Y$ and $R^ih_*\mathcal{O}_Z=0$ for all $i>0$. Hence $h$ is an rational resolution of $Y$, so $Y$ also has rational singularity.

\end{proof}

Next we proof that flips also preserves rational singularity.

To proof this theorem first we will need a lemma 
\begin{lemma}\label{T2}
Let $(X,\Delta)$ be a lc pair and let $g:(Z,\Delta_Z)\rightarrow (X,\Delta)$ a dlt model of $(X,\Delta)$. Let $f:(Y,\Delta_Y)\rightarrow (X,\Delta)$ be a log resolution of $(X,\Delta)$, which is also a log resolution of $(Z,\Delta_Z)$. Then every associated prime of $R^if_*\mathcal{O}_Y$ is the generic point of some  log canonical center of $(X,\Delta)$ for every $i>0$	
	
\end{lemma}

\begin{proof}
We take a dlt model $g:(Z,\Delta_Z)\rightarrow (X,\Delta)$ which is possible by \cite[Lemma $2.23$]{DH22} 
So we have that $K_Z+\Delta_Z=g^*(K_X+\Delta)$ and $(Z,\Delta_Z)$ is a K\"ahler $3$-fold dlt pair. We take projective bimeromorphic morphism $h:Y\rightarrow Z$ such that $K_Y+\Delta_Y=h^*(K_Z+\Delta_Z)$, with $Y$ smooth and $\Supp \Delta_Y$ simple normal crossing divisor on $Y$. We may assume that log canonical center of $(Z,\Delta_Z)$ is contained in the stratum. Then we have \[
	K_Y+\Delta_Y\sim_{\mbR,f} 0.
\]
 
 Where  $f=g\circ h:Y\rightarrow X$.
 
$\Delta_Y$ is a SNC divisor. We will apply \cite[theorem 1.1(i)]{FV22}. Note that $-(K_Y+\Delta_Y)$ is $f$-trivial hence $f$-semiample. Hence by the above mentioned theorem we have that every associated prime of $R^if_*\mathcal{O}_Y$ is supported on $f$ image of strata of $(Y,\Delta)$. Now all the lc centeres of $(Y,\Delta_Y)$ are contained in the strata, whose images are the generic points point of log canonical centers of $(X,\Delta)$. Hence every associated prime of $R^if_*\mathcal{O}_Y$ is the generic point of some log canonical center of $(X,\Delta)$ for every $i>0$.

\end{proof}
\begin{proof}[Proof of Theorem \ref{T4}](2)
	Let $g:Z\rightarrow X^+$ be a resolution as in Lemma \ref{T2}. Let $\Ex(f^+)$ be the exceptional locus of $f^+:X^+\rightarrow Y$. As $X\setminus \Ex(f)$ is isomorphic to $X^+\setminus \Ex(f^+)$ hence $X^+\setminus \Ex(f^+)$ has only rational singularities. Thus $\Supp R^ig_* \mathcal{O}_Z\subset \Ex(f^+)$ for every $i>0$ since outside $\Ex(f^+)$ these sheaves are $0$. By negativity lemma there are no log canonical centeres of $(X^+,\Delta^+)$ contained in $\Ex(f^+)$. By Lemma \ref{T2} every associated prime of $R^ig_*\mathcal{O}_Z$ is contained in generic point of some log canonical center of $(X^+,\Delta^+)$ for every $i>0$. Thus we have $R^ig_*\mathcal{O}_Z=0$ for every $i>0$. Thus $g$ is a rational resolution. Hence $X^+$ has only rational singularities by \cite[Lemma $5.12$]{KM}.
	
\end{proof}

\section{Termination}
Now we come to proving termination of flips for log canonical K\"ahler $3$-fold pairs.

\begin{proof}[Proof of Theorem \ref{T1}]
Similar proof as \cite[Theorem $8.2$]{AST} for algebraic $3$-folds works here. Let us emphasize the proof a bit more.\\
Let 
\begin{equation*} 
\xymatrix{     
               (X_0,B_0)\ar[dr]_{\phi_0} \ar@{-->}[rr] &     & (X_1,B_1)\ar[dl]^{\phi_0+}\ar[dr]_{\phi_1} \ar@{-->}[rr]& & (X_2,B_2)&\\
                                             &Z_0&                  & Z_1&
                                             }
 \end{equation*} 	
	be a sequence of flips let $q_0:(Y_0,D_0)\rightarrow (X_0,B_0)$ be a $\mbQ$-factorial dlt model, so $K_{Y0}+D_0=q_0^*(K_{X0}+B_0)$. So $K_{Y0}+D_0$ is not nef over $Z_0$. Then by \cite[Theorem $1.1$]{DH22} there is a sequence of divisorial contractions and flips which terminates, and at the end we get $(Y_1,D_1)\rightarrow Z_0$ a dlt pair with $K_{Y1}+D_1$ is relatively nef. Note that this is true even if the starting $K_{Y_0}+D_0$ is not pseudo-effective as we are running relative MMP for birational morphisms, $K_{Y_0}+D_0$ is always relatively pseudo-effective. So by definition $(Y_1,D_1)\rightarrow Z_0$ is a weak canonical model of $(Y_0,D_0)\rightarrow Z_0$ in the sense of \cite[Definition $3.6.7$]{BCHM}. The contracted locus by MMP $Y_0\dashrightarrow Y_1$ is contained in $G\cup E$, where image of $E$ in $X_0$ is the flipping locus, and $G$ is the support of the exceptional locus of $Y_0\rightarrow X_0$. If $C\subset E$ is a curve vertical over $Z_0$ then $(K_{Y_0}+D_0)\cdot C<0$. As $K_{Y_1}+D_1$ is nef, we have the whole locus $E$ is contracted by this MMP. So the exceptional locus of $(Y_1,D_1)\rightarrow Z_0$ is contained in image of $G$ in $Y_1$. Note that as $G\subset \supp D_0$ is union of all the lc center of $(Y_0,D_0)$, union of all the lc centers of $(Y_1,D_1)$ is the image of $G\subset \supp D_1$. Now let $F$ be an surface getting contracted to a point by the morphism $Y_1\rightarrow Z_0$. Then $F$ is contained in the image of $G$, hence $F$ is a lc centere of $(Y_1,D_1)$ and contained in $\Supp{D_1}$. By adjunction we have that $(K_{Y_1}+D_1)|_F=K_F+\theta$ is also nef. By \cite[Theorem $6.2$]{DO22}, $K_F+\theta$ is semi-ample. Now suppose $C$ is a curve getting contracted by $Y_1\rightarrow Z_0$, then $C$ is also an lower dimensional lc center as it is contained in image of $G$, hence contained in the strata of $D_1$. Hence there exists a lc center $S$ of dimension $2$ which contains $C$. Again by adjunction and abundance for compact K\"ahler lc surfaces we get $(K_{Y_1}+D_1)|_S=K_S+\theta'$ is semi-ample, thus $(K_S+\theta')|_C$ is also semi-ample. So $K_{Y_1}+D_1$ is relatively semi-ample over $Z_0$. We have the morphism $q_1:(Y_1,D_1)\rightarrow (X_1,B_1) $ by uniqueness of log canonical model over $Z_0$. So we have that $q_1^*(K_{X_1}+B_1)=K_{Y_1}+D_1$, hence $(Y_1,D_1)$ is also the dlt model of $(X_1,B_1)$. Now we repeat the same process on the map $Y_1\rightarrow X_1\rightarrow Z_1$. This way a sequence of flips on $X$ lifts to a sequence of flips and divisorial contraction on $Y_0$. As $(Y_0,D_0)$ is compact K\"ahler $3$-fold dlt pair, this sequence terminates as this is also a $K_{Y_0}+D_{0}$-MMP, and hence the sequence of flips below is also finite.
\end{proof}

\section{MMP for Uniruled Pairs}
In this section we consider Minimal Model Program for non-pseudo-effective lc compact K\"ahler $3$-fold pairs $(X,B)$. Since $K_X+B$ is not pseudo-effective, neither $K_X$ is, the MRC fibration $X\dashrightarrow Z$ is non trivial. If $\dim Z<2$, Then by \cite[Lemma $2.39$]{DH22} $X$ is projective. Since the projective case is well understood we focus on the case where $\dim Z=2$. Also note that $Z$ is not uniruled hence $K_Z$ is pseudo-effective. Moreover $Z$ is in the Fujiki class $\mathcal{C}$. Then replacing $Z$ by a resolution of singularities we may assume that $Z$ is a smooth compact complex surface in Fujiki's class $\mathcal{C}$, hence $Z$ is also K\"ahler. 
\begin{definition}
Let $(X,B\geqslant 0)$ be a log pair, where $X$ is $\mbQ$-factorial compact K\"ahler $3$-fold. Suppose that the base of the MRC fibration $f:X\dashrightarrow Z$ has dimension $2$. Let $X_z\cong \mbP^1$ be a general fiber of $f$. Then a modified K\"ahler class $\omega$ on $X$ is called $K_X+B$ normalized if $(K_X+B+\omega).X_z=0$. Note that $\omega$ is modified K\"ahler class, it is positive on general fibers $X_z$, and hence $(K_X+B)\cdot X_z<0$	
	
\end{definition}

\begin{lemma}\label{l1}
With the same hypothesis as in the defintion above assume that $(X,B)$ is lc and $\omega$ is a $K_X+B$ normalised modified K\"ahler class. Then $K_X+B+\omega$ is pseudoeffective.	
	
\end{lemma}
\begin{proof}
The same proof of \cite[Lemma $4.2$]{DH22} works in this case since the only fact about dlt which is used in the proof is that discrepency of any exceptional divisor over $X$ is greater than equals to $-1$, which is also true in lc singularities.	
	
\end{proof}

Next we prove Cone theorem in uniruled set up.
\begin{theorem}\label{C2}
	Let $(X,B)$ be a lc pair, where $X$ is a $\mbQ$-factorial compact K\"ahler $3$-fold with rational singularities. Suppose that $X$ is uniruled and the base of MRC fibration $X\dashrightarrow Z $ is a surface. Suppose $\omega$ be a modified $K_X+B$ normalised K\"ahler class. Then there exists a countable family of curves $\Gamma_i$ on $X$ and a positive number $d$ such that $0<-(K_X+B+\omega)\cdot \Gamma_i\leqslant d$ and \[
		\overline{NA}(X)=\overline{NA}(X)_{(K_X+B+\omega)\geqslant 0}+\sum_{i\in I}\mbR^+[\Gamma_i]
	\]
\end{theorem}
\begin{proof}
	By \cite[Theorem $4.6$]{DH22} we know the cone theorem for dlt pairs in the same set-up. So let us go to a log resolution and running MMP for dlt pair we get to the dlt model of $(X,B)$, let that be $(X',B')$. Now note that if $\omega$ is a modified K\"ahler class $(K_X+B)$-normalised, then pullback of $\omega$ is $K_X'+B'$ normalised big class. Then using same technique as Theorem \ref{C1} we get our desired cone theorem. 
	
\end{proof}

Now we get to the contraction theorems.
\begin{theorem}\label{C'2}
Let $(X,B)$ be a lc pair, where $X$ is a $\mbQ$-factorial compact K\"ahler $3$-fold, with klt singularity. Suppose that $X$ is uniruled and the dimension of base of the MRC fibration $X\dashrightarrow Z$ is 2, and $(K_X+B)\cdot F<0$ for a general fiber $F$ of $X\dashrightarrow Z$. Let $\omega$ be a Kahler class such that $K_X+B+\omega$ is pseudo-effective, and $R$ is a $(K_X+B+\omega)$-negative extremal ray $R$. Then the contraction $c_R:X\rightarrow Y$ of the ray $R$ exists.
	
\end{theorem}
\begin{proof}
We may assume that the extremal ray $R$ is cut out by a $(1,1)$ nef class $\alpha$. Rescaling $\alpha$ if necessary, we see that $\eta=\alpha-(K_X+B+\omega)$ is positive on $\overline{NA}(X)\setminus 0$ by Theorem \ref{C2}. Thus it follows that $\eta$ is K\"ahler class, so $\alpha$ is a nef and big class, since $K_X+B+\omega$ is pseudo-effective.\\
Suppose $R$ is small. Then exact same argument for the uniruled portion of Theorem \ref{C'1} proves the contraction.\\
Suppose that $R$ is of divisorial type, then the corresponding irreducible divisor $S$ is covered by and contains all the curves $C\subset X$ such that $[C]\in R$. Now if nef dimension of $\alpha=0$, where $\alpha$ is the class whose null locus is exactly $S$, then we have the contraction just like the pseudo-effective case, even without assuming $X$ having klt singularity. But for the case nef dimension of $\alpha=1$, we will need to assume $X$ has klt singularity. The proof is exactly same as Theorem \ref{Cont}.

\end{proof}

\section{Existence of Mori fiber space}
In this section we prove the existence of mori fiber space for $\mbQ$-factorial log canonical uniruled pairs. We will need some theorems first.
\begin{theorem}\label{NPM1}
Let $(X,B)$ be a lc pair, where $X$ is a $\mbQ$ factorial compact K\"ahler $3$-fold. Suppose that $X$ is uniruled and the base of the MRC fibration $f:X\dashrightarrow Z$ has dimension 2 and $(K_X+B)\cdot F<0$, where $F$ is a general fiber of $f$. Then there is a bimeromorphic map $\phi:X\dashrightarrow X'$ given by a sequence of $K_X+B$-flips and divisorial contractions such that for any $(K_{X'}+B')$-normalised K\"ahler class $\omega'$ on $X'$, the adjoint class $K_{X'}+B'+\omega'$ is nef. Where $B'=\phi_*B$.
\end{theorem}
\begin{proof}
	Suppose that there is a $(K_X+B)$-normalised K\"ahler class $\omega$ on $X$, such that the adjoint class $K_X+B+\omega$ is not nef. But by Lemma \ref{l1} we know that $K_X+B+\omega$ is pseudo effective. Then applying Theorem \ref{C'2} we know that we can contract $R=(K_X+B+\omega)^{\bot}\cap \overline{NA}(X)$. Note that $R$ is also $K_X+B$ negative. Hence applying Theorem \ref{T1} we see that this MMP terminates with $(X',B')$ such that $(K_{X'}+B'+\omega')$ is nef.
	
\end{proof}

Now we will state the theorem which is going to be the crux arguement for the existence of Mori fiber space.
\begin{proposition}\cite[Corollary $5.4$]{DH22}\label{MC1}
Let $(Y,B_Y)$ be a klt pair, where $Y$ is a $\mbQ$-factorial compact K\"ahler $3$-fold. Suppose that $Y$ is uniruled and the base of the MRC fibration $g:Y\dashrightarrow Z'$ has dimension $2$. Let $\omega_Y$ be a nef and big class on $Y$ such that $K_Y+B_Y+\omega_Y$ is nef and $(K_Y+B_Y+\omega_Y)\cdot F=0$, where $F$ is a general fiber of $g$. Then there exists a proper surjective map with connected fibers $\psi:Y\rightarrow S$ onto a compact K\"ahler normal surface $S$ such that $K_Y+B_Y+\omega_Y$ is $\psi$-trivial.  	
	
\end{proposition}
	
Now we come to proof of our main theorem of this section

\begin{proof}[Proof of Theorem \ref{FT}]
	 Since $K_X+B$ is not pseudo effective, $K_X$ is not pseudo-effective. Thus by Brunellas theorem applied to a resolution of $X$ it follows that $X$ is uniruled. We may assume that the dimension of the base of the MRC fibration $f:X\dashrightarrow Z$ is 2, as otherwise $X$ would be projective. Let $F$ be a general fiber of the MRC fibation $f$. By Lemma \ref{l1} if $(K_X+B)\cdot F\geqslant 0$, then $K_X+B$ is pseudo-effective, contradicting our assumptions. Therefore $(K_X+B)\cdot F< 0$.  Then by Theorem \ref{NPM1} there is a $(K_X+B)$-MMP, $\psi: X\dashrightarrow X'$ such that for every $(K_{X'}+B')$-normalized K\"ahler class $\omega'$, the class $K_{X'}+B'+\omega'$ is nef, where $B'=\psi_* B$. Since $X'$ is K\"ahler and $(K_{X'}+B')\cdot F'<0$, where $F'\cong F$ is a general fiber of the induced MRC fibration $X'\dashrightarrow Z$, we may pick a $(K_{X'}+B')$-normalized K\"ahler class, say $\omega'$. We have that $K_{X'}+B'+\omega'=K_{X'}+(1-\epsilon)B'+\omega+\epsilon B'$, where $\epsilon>0$ is small enough so that $\omega+\epsilon B'$ is K\"ahler. Now as $X'$ is klt, so is $(X',(1-\epsilon)B')$ By proposition \ref{MC1}, there exists holomorphic fibration $\Psi:X'\dashrightarrow S'$ on to a normal compact K\"ahler surface $S'$ such that $K_{X'}+B'+\omega'$ is $\Psi$-trivial. In particular $\Psi$ is a projective morphism hence the theorem follows for usual relative minimal model program for projective morphisms.
	
\end{proof}
Now we prove the final theorem of this article which is existence and termination of MMP for compact K\"ahler $3$-fold lc pair $(X,\Delta)$, with $X$ having klt singularity.
\begin{proof}[Proof of \ref{MMP}]
$(1)$ By the Theorems \ref{C1},\ref{Cm},\ref{T4}, \ref{T1} we have that for $K_X+B$ pseudo-effective any sequence of $K_X+B$ divisorial contraction and flips terminates.\\
$(2)$ From proof of Theorem \ref{FT} and the statement of Theorem \ref{NPM1} notice that all the MMP we are runnning to reach a Mori fiber space is in-fact $K_X+\Delta$-MMP. And by Theorem \ref{T1} we have that any such MMP will terminate. As we have already shown existence of Mori fiber space, we have any such $K_X+\Delta$-MMP with $K_X+\Delta$ not pseudo-effective will terminate with a Mori fiber space.
    
\end{proof}

\bibliographystyle{hep}
\bibliography{bibliography}

\end{document}